\documentclass[10pt]{amsart}
\input{epsf}
\usepackage{amssymb,latexsym}
\topskip  \numberwithin{equation}{section}
\newtheorem{theorem}{Theorem}[section]

\newtheorem{corollary}[theorem]{Corollary}

\newtheorem{lemma}[theorem]{Lemma}

\newtheorem{example}[theorem]{Example}

\begin{document}

\pagenumbering{arabic}
\pagestyle{headings}
\def\sof{\hfill\rule{2mm}{2mm}}
\def\llim{\lim_{n\rightarrow\infty}}

\title{Partially ordered patterns and compositions}
\maketitle

\begin{center}
{\bf Silvia Heubach}\\
{\it Dept. of Mathematics, California State University Los Angeles,
Los Angeles, CA 90032, USA}\\
{\tt sheubac@calstatela.edu}\\
\vskip 10pt
{\bf Sergey Kitaev}\\
{\it Department of Mathematics, Reykjavik University,
IS-103 Reykjavik, Iceland}\\
{\tt sergey@ru.is}\\
\vskip 10pt
{\bf Toufik Mansour}\\
{\it Department of Mathematics, Haifa University, 31905 Haifa, Israel}\\
{\tt toufik@math.haifa.ac.il}\\

\end{center}

\def\P{POP}
\def\A{\mathcal{A}}
\def\SS{\frak S}
\def\NN{\mathbb N}
\def\Ps{POPs}
\def\mn{\mbox{-}}
\def\newop#1{\expandafter\def\csname #1\endcsname{\mathop{\rm #1}\nolimits}}
\newop{MND}
\section*{Abstract}
A partially ordered (generalized) pattern (POP) is a  generalized pattern some of whose letters are incomparable, an extension of generalized permutation patterns introduced by
Babson and Steingr{\'\i}msson. POPs were introduced in the symmetric group by Kitaev~\cite{Ki2, Ki4}, and studied in the set of $k$-ary words by Kitaev and Mansour~\cite{KM}. Moreover, Kitaev et al.~\cite{KMP} introduced segmented POPs  in compositions. In this paper,
we study avoidance of POPs  in compositions and generalize results for avoidance of POPs in permutations and words. Specifically, we obtain results for the generating functions for the number of compositions that avoid shuffle patterns and multi-patterns. In addition, we give the
generating function for the distribution of the maximum number of
non-overlapping occurrences of a segmented POP $\tau$ (that is
allowed to have repeated letters) among the compositions of $n$
with $m$ parts in a given set, provided we know the generating function for the
number of compositions of $n$ with $m$ parts in the given set that
avoid $\tau$. This result is a $q$-analogue of the main result
in~\cite{KM}.

\noindent{\bf Keywords}: Compositions, partially ordered
(generalized) patterns,  non-overlapping
occurrences, generating functions.

\noindent{\bf 2000 Mathematics Subject Classification}: 05A05, 05A15
\thispagestyle{empty}
\section{Introduction}

Pattern avoidance was originally studied in permutations
(see~\cite{bona, SS}), and the patterns studied were also
permutation patterns. Generalizations in several directions took
place: 1) Looking at pattern avoidance in permutations with
different types of patterns and avoidance of sets of patterns (see
\cite{ Ki2} and references therein), and 2) asking the same
questions for words (see \cite{Burstein, BM2, BM3, KM}).
Independently, several authors (see \cite{CGH, CH, CH2, G, G2, HB,
HM, HM2}) gave results on enumerating compositions of $n$ with parts
in a given set $A$ according to rises, levels and drops (which can
be considered as the simplest 2-letter patterns). Heubach and
Mansour (see \cite{HM3, HM4}) combined these two areas by giving
results on the generating function for the number of compositions of
$n$ with $m$ parts in a set $A$ that  avoid 3-letter patterns.
Moreover, Kitaev et al.~\cite{KMP} introduced segmented partially
ordered (generalized) patterns in compositions.

In this paper we generalize some of the results in the literature on pattern avoidance in
permutations~\cite{Ki2}, \cite{Ki4} and words~\cite{KM} by
studying pattern avoidance of partially ordered patterns {\rm(}\Ps{\rm)} in
compositions. Section~\ref{basics} contains basic definitions and terminology. In Section~\ref{main}, we give a
general result that expresses the generating function of the number of compositions that avoid a POP composed of two smaller patterns in terms of the generating functions for the smaller patterns. We apply this result (Theorem~\ref{thmpc}) to two specific types of POPs, namely shuffle patterns and multi-patterns and show equivalence for families of patterns of each type. We close in Section~\ref{non-overlap} by giving a result for the maximum number of non-overlapping occurrences of a \P\, in a composition, which is a generalization of a theorem proved by Kitaev  ~\cite[Theorem 32]{Ki2} for permutations and by Kitaev and Mansour \cite[Theorem 5.1]{KM} for words.

\section{Definitions and Terminology} \label{basics}

Let $\NN$ be the set of all positive integers, and let $A$ be any ordered finite set
of positive integers, say $A=\{a_1,a_2,\ldots,a_k\}$, where
$a_1<a_2<a_3<\cdots<a_k$. (An ``ordered set'' in this paper will
always refer to a set whose elements are listed in increasing order.) Also, let $[k]^n$ denote the set of all words of length $n$ over the
(totally ordered) alphabet $[k]=\{ 1,2,\dots,k \}$.

A {\em composition} $\sigma=\sigma_1\sigma_2\ldots\sigma_m$ of
$n\in\NN$ is an ordered collection of one or more positive integers
whose sum is $n$. The number of {\em summands} or {\em letters},
namely $m$, is called the number of {\em parts} of the composition.
For any ordered set $A=\{a_1,a_2,\ldots,a_k\}\subseteq\NN$, we
denote the set of all compositions of $n$ with parts in $A$ (resp.
with $m$ parts in $A$) by $C_n^A$ (resp. $C_{n;m}^A$).

 A \emph{generalized pattern}  $\tau$ is a word in $[\ell]^m$ (possibly with dashes between
some letters) that contains each letter from $[\ell]$ (possibly with
repetitions). Generalized patterns that contain dashes in
all possible positions (e.g., $2\mn1\mn4\mn3$) are called {\em classical patterns}. Note that classical patterns place no adjacency requirements on occurrences of the letters of
a pattern in words or compositions.  If all the dashes are removed, we
have a {\em consecutive}, or {\em segmented}, {\em pattern}. For ease of readability, we will refer to generalized patterns simply as patterns in the remainder of this paper.

We say that a composition
$\sigma\in C_n^A$ \emph{contains} a pattern $\tau$ if $\sigma$ contains
a subsequence isomorphic to $\tau$ in which the entries
corresponding to consecutive entries of $\tau$ (those not
separated by a dash) must be adjacent. Otherwise, we say that
$\sigma$ \emph{avoids} $\tau$ and write $\sigma\in C_n^A(\tau)$. Thus,
$C_n^A(\tau)$ denotes the set of all
compositions of $n$ with parts in $A$ that avoid $\tau$. Moreover,
if $T$ is a set of patterns, then
$C_n^A(T)$ denotes the set of all
compositions of $n$ with parts in $A$ that avoid each pattern from
$T$ simultaneously. For example, $241874$ avoids $312$ and contains
three occurrences of $1\mn32$, namely $287$, $274$ and $487$. (Note
that 284 is not an occurrence of $1\mn32$ due to the adjacency
requirement).

Kitaev~\cite{Ki2},~\cite{Ki4} introduced {\em partially ordered patterns} ({\em POPs})\footnote{In ~\cite{Ki2},
POPs are called {\em POGPs} ({\em Partially Ordered Generalized
Patterns}). We use POPs instead to shorten the notation in this paper.}  on permutations, which
extend  {\em generalized permutation patterns}  introduced by Babson and Steingr{\'\i}msson~\cite{BS}. Specifically, a \P\ $\tau$ is a word consisting of letters from a
partially ordered
    alphabet $\mathcal{T}$ such that the letters in $\tau$
    constitute an order ideal in $\mathcal{T}$. If letters $a$ and
$b$ are incomparable in a POP $\tau$, then the relative size of the letters in
$\sigma$ corresponding to $a$ and $b$ is unimportant in an occurrence of $\tau$
in  $\sigma$. For instance,
if $\mathcal{T} = \{ 1, 1', 2'\}$ and the only relation is $1' < 2'$, then the sequence
    31254 has two occurrences of $\tau =1 1' 2'$, namely $3 1 2$ and $1 2
    5$. As for  generalized patterns, if a \P\ $\tau=\tau_1\ldots \tau_k$ has a dash between, say, $\tau_i$
    and  $\tau_{i+1}$, then in an occurrence of $\tau$ in a
composition $\sigma$, the letters corresponding to $\tau_i$
    and  $\tau_{i+1}$ do not have to be adjacent. For example, for $\mathcal{T}$ given
    above, if $\tau=1\mn1'2'$, then the composition $113425$
 contains seven occurrences of $\tau$,
namely $113$, $134$ twice, $125$ twice, $325$, and $425$.

Following~\cite{Ki2} and ~\cite{KM}, we consider two particular
classes of \Ps\, --  {shuffle patterns} and {multi-patterns} --
which allows us to give an analogue of the main results
in~\cite{Ki2} and ~\cite{KM} for compositions. Let $\{\tau_0,\tau_1,\dots,\tau_s\}$ be a set of consecutive
patterns. A {\em multi-pattern} is of
the form $\tau=\tau_1\mn\tau_2\mn\cdots\mn\tau_s$ and a {\em shuffle
pattern} is of the form $\tau=\tau_0\mn a_1\mn\tau_1\mn
a_2\mn\cdots\mn\tau_{s-1}\mn a_s\mn\tau_s$, where  each
letter of $\tau_i$ is incomparable with any letter of $\tau_j$ whenever $i\neq j$. In addition, the letters $a_i$ are either
all greater or all smaller than any letter of $\tau_j$ for any $i$ and $ j$. For example, $1'\mn2\mn1''$ is a shuffle pattern, and $1'\mn1''$ is
a multi-pattern. Clearly, we can get a multi-pattern from a
shuffle pattern by removing all the letters $a_i$. Furthermore, there is a connection
between avoidance of a \P\ and multi-avoidance of generalized patterns in compositions. For example, avoiding the \P\, $2'\mn1\mn2''$ is the same as simultaneously avoiding
 the patterns $2\mn1\mn2$, $3\mn1\mn2$, and
$2\mn1\mn3$~(similar to \cite[Proposition~2.7]{KM}).

\section{POPs in compositions with parts in a given set}\label{main}
We will now derive results on avoidance of \Ps\, in compositions. In order to distinguish which letters are comparable and which ones are not, we will use primes in the following way. If two
letters, say 1 and 2, have the same number of primes, say two, then
they are comparable and naturally $1''<2''$. Any two letters with a
different number of primes are incomparable. Unless dealing with
shuffle or multi-patterns, if a letter in a POP has no primes, then that
 letter is greater than every letter with one or more
primes and we will emphasize this fact by
using a value that is bigger than those for the primed letters.
For example, in $\tau=1'\mn 2\mn 1''$, the second letter is the
greatest one and the first and the last letters are incomparable to
each other. The composition $\sigma=31421$
has five occurrences of $\tau$, namely $342$, $341$, $142$, $141$,
and $121$.

Let $C_{\tau}^A(x)=\sum\limits_{n\geq0}|C_n^A(\tau)|x^n$ (resp.
$C_\tau^A(x;m)=\sum\limits_{n\geq0} |C_{n;m}^A(\tau)|x^n$ and
$C_{\tau}^A(x,y)=\sum\limits_{n,m\geq0}|C_{n;m}^A(\tau)|x^ny^m$)
denote the generating function for the numbers $|C_n^A(\tau)|$ (resp.
$|C_{n;m}^A(\tau)|$) of compositions in $C_n^A$ (resp. $C_{n;m}^A$)
avoiding the pattern $\tau$. For example, if $A=\{a_1,a_2,\ldots,a_k\}$ is any ordered set and
$\tau=1'\mn2\mn1''$, then we have
\begin{equation}
C_{1'\mn2\mn1''}^A(x,y) = \frac{1}{\prod_{a\in A}(1-x^ay)^2}
-\sum_{a\in A}\frac{x^ay}{\prod_{a\leq b\in
A}(1-x^by)^2}.\label{eqex1}
\end{equation}
This result follows from the specific structure of the compositions $\sigma$ that avoid $\tau=1'\mn2\mn1''$. If $\sigma$ avoids $\tau$, and $\sigma$
contains $s>0$ copies of the letter $a_k$, then the letters $a_k$ can only appear as blocks on the left and right end
 of $\sigma$. If $\sigma$ contains no $a_k$, then
$\sigma\in C_{n,m}^{A'}(\tau)$ where $A'=A-\{a_k\}$. So, for all
$n\geq0$, we have
$$C_{\tau}^A(x;m)=\sum_{i=0}^{m-1}(i+1)x^{ia_k}C_\tau^{A'}(x;m-i)+x^{ma_k},$$
since the generating function for the possibilities to place $i$ letters $a_k$ into
$\sigma$ is given by $(i+1)C_\tau^{A'}(x;m-i)$, for $0\leq i <
m$, and by $x^{m\,a_k}$ for $i = m$. Thus, for $m \geq 2$,
$$
C_\tau^A(x;m)-2x^{a_k}C_\tau^A(x;m-1)
=C_\tau^{A'}(x;m)-x^{2a_k}\left(\sum_{s=0}^{m-3}(s+1)x^{sa_k}C_\tau^{A'}(x;m-2-s)-x^{(m-2)a_k}\right)$$
or equivalently,
$$
C_\tau^A(x;m)-2x^{a_k}C_\tau^A(x;m-1)+x^{2a_k}C_\tau^{A}(x;m-2)=C_\tau^{A'}(x;m),$$

together with $C_\tau^A(x;0)=1$ and $C_\tau^A(x;1)=\sum_{a\in
A}x^ay$. Multiplying both sides of the recurrence above by $y^m$,
summing over all $m\geq2$ and using induction on elements of $A$
together with the fact that
$C_\tau^{\{a_1\}}(x,y)=\frac{1}{1-x^{a_1}y}=\frac{1}{(1-x^{a_1}y)^2}-\frac{x^{a_1}y}{(1-x^{a_1}y)^2}$,
we get (\ref{eqex1}). Equation~(\ref{eqex1}) for $x=1$ and $A=[k]$
gives the corresponding result for words ~\cite[Equation~2.1]{KM}.

In order to prove general results, it is  convenient to introduce the
notion of quasi-avoidance. Let $\tau$ be a consecutive pattern. A composition $\sigma$ {\em
quasi-avoids} $\tau$ if $\sigma$ has exactly one occurrence of
$\tau$ and this occurrence consists of the $|\tau|$ rightmost parts
of $\sigma$, where $|\tau|$ denotes the number of letters in $\tau$. For example, the composition $4112234$
quasi-avoids the pattern $1123$, whereas the compositions
$5223411$ and $1123346$ do not.

First,  relate the generating function for the number of compositions avoiding a given pattern $\tau$ with the generating function for the number of compositions that quasi-avoid $\tau$.

\begin{lemma}\label{praa}
Let $\tau$ be a non-empty consecutive pattern. Let $D_\tau^{A}(x,y)$
denote the generating function for the number of compositions in $C_{n;m}^A$
that quasi-avoid $\tau$. Then
\begin{equation}
D_\tau^{A}(x,y)=1+C_\tau^A(x,y)\left(y\sum_{a\in A}x^a
\,-1\right).
\end{equation}
\end{lemma}
\begin{proof}
We use  arguments similar to  in the proof of
\cite[Proposition~4]{Ki2}. Adding the part $a$ to a composition with $m-1$ parts that avoids $\tau$ creates either a composition with $m$ parts that still avoids $\tau$ or that quasi-avoids $\tau$.  Thus, for $m\geq 1$,
$$D_{\tau}^{A}(x;m)=\left(\sum_{a\in A}x^a\right) C_\tau^A(x;m-1)-C_\tau^A(x;m).$$
 Multiplying both sides of this equality by $y^m$ and summing
over all natural numbers $m$ we get the desired result.
\end{proof}

Lemma~\ref{praa} for $A=[k]$ and $x=1$ gives the corresponding result for words
\cite[Proposition~2.4]{KM}.

We now obtain a general theorem that is a good auxiliary tool for calculating the
generating function for the number of compositions that avoid a given \P.

\begin{theorem}\label{thmpc}
Let $A=\{a_1,\ldots,a_k\}$ be any ordered finite set of positive
integers. Suppose $\tau=\tau_0\mn\phi$, where $\phi$ is an
arbitrary \P, and the letters of $\tau_0$ are incomparable to the
letters of $\phi$. Then for all $k\geq 1$, we have
$$C_{\tau}^A(x,y)=C_{\tau_0}^A(x,y)+C_\phi^A(x,y)D_{\tau_0}^{A}(x,y).$$
\end{theorem}
\begin{proof}
 To find $C_\tau^A(x,y)$, we observe
that there are two possibilities: either $\sigma$ avoids $\tau_0$,
or $\sigma$ does not avoid $\tau_0$. In the first case, the generating function is given by $C_{\tau_0}^A(x,y)$. If $\sigma$ does not avoid $\tau_0$, then we can write $\sigma$ in the form $\sigma=\sigma_1\sigma_2\sigma_3$,
 where $\sigma_1\sigma_2$ quasi-avoids the pattern
$\tau_0$, and $\sigma_2$ is order isomorphic to $\tau_0$. Clearly,
$\sigma_3$ must avoid $\phi$, thus, the generating function is equal to
$C_\phi^A(x,y)D_{\tau_0}^{A}(x,y)$, and we obtain the stated result.
\end{proof}

Theorem~\ref{thmpc} can be used for reduction, but also to easily compute the generating function for a new pattern from the generating function of a known pattern. We will use the notion of equivalence of patterns to obtain several results that hold for whole families of patterns. Two \Ps  \,\,$\tau$ and $\phi$ are said to be {\em equivalent}, and we write $\tau\equiv\phi$, if the number of compositions in $C_{n;m}^A$  that
avoid $\tau$ is equal to the number of compositions in $C_{n;m}^A$ that
avoid $\phi$ for all $n,m$.

Let $A=\{a_1,\ldots,a_k\}$ be any ordered finite set of positive
integers and $\sigma=\sigma_1\sigma_2\ldots\sigma_m\in C_{n;m}^A$. Then
the {\em reverse $R(\sigma)$} of a composition $\sigma$ is the
composition $\sigma_m\ldots\sigma_2\sigma_1$.
We call this bijection of $C_{n;m}^A$ to itself {\em
trivial}. (The other trivial bijection is $I$, the identity bijection).  Note that the complement operation defined for permutations and words is not defined for
compositions. It is easy to see that $\tau\equiv R(\tau)$ for any pattern $\tau$.
For example, the number of compositions that avoid the pattern
$21\mn2$ is the same as the number of compositions that avoid the pattern
$2\mn12$.

In the following two subsections we obtain results for two specific classes of \Ps\,  -- shuffle patterns and multi-patterns.

\subsection{Shuffle patterns in compositions.}\label{sec3}
We  consider the shuffle patterns $\tau\mn \ell\mn \nu$ and  $\tau\mn1\mn \nu$, where
$\ell$ (resp. $1$) is the greatest (resp. smallest) element of the
pattern.

\begin{theorem}\label{thsb}
Let $A=\{a_1,a_2,\ldots,a_k\}$ be any ordered set of positive
integers.

\begin{enumerate}
\item Let $\phi$ be the shuffle pattern $\tau\mn \ell \mn\nu$. Then
for all $k\geq \ell$,
$$C_\phi^A(x,y)=\frac{C_\phi^{A-\{a_k\}}(x,y)-x^{a_k}yC_\tau^{A-\{a_k\}}(x,y)C_\nu^{A-\{a_k\}}(x,y)}
{(1-x^{a_k}yC_\tau^{A-\{a_k\}}(x,y))(1-x^{a_k}yC_\nu^{A-\{a_k\}}(x,y))}.$$

\item Let $\psi$ be the shuffle pattern $\tau\mn1\mn\nu$. Then
for all $k\geq \ell$,
$$C_\psi^A(x,y)=\frac{C_\psi^{A-\{a_1\}}(x,y)-x^{a_1}yC_\tau^{A-\{a_1\}}(x,y)C_\nu^{A-\{a_1\}}(x,y)}
{(1-x^{a_1}yC_\tau^{A-\{a_1\}}(x,y))(1-x^{a_1}yC_\nu^{A-\{a_1\}}(x,y))}.$$
\end{enumerate}
\end{theorem}

\begin{proof}
We derive a recurrence relation for $C_\phi^A(x,y)$  where $\phi=\tau\mn \ell \mn\nu$.
Let $\sigma\in C_{n,m}^A(\phi)$ be
such that it contains exactly $s$ copies of the letter $a_k$. If
$s=0$, then the generating function for the number of such compositions is
$C_\phi^{A'}(x,y)$, where $A'=A-\{a_k\}$. For $s \geq 1$, we write
$\sigma=\sigma_0 a_k\sigma_1 a_k\cdots a_k\sigma_{s}$,
where $\sigma_j$ is a $\phi$-avoiding composition with parts in
$A'$, for $j=0,1, \ldots, s$. Then either
$\sigma_j$ avoids $\tau$ for all $j$, or there exists a $j_0$ such
that $\sigma_{j_0}$ contains $\tau$, $\sigma_j$ avoids $\tau$ for
all $j=0,1, \ldots, j_0-1$ and $\sigma_j$ avoids~$\nu$ for any
$j=j_0+1,\ldots,s$. In the first case, the generating function for the number of such
compositions is
$x^{sa_k}y^s\left(C_\tau^{A'}(x,y)\right)^{s+1}$. In the second
case, the generating function is given by
$$x^{sa_k}y^s\sum_{j=0}^s\left(C_\tau^{A'}(x,y)\right)^j\left(C_\nu^{A'}(x,y)\right)^{s-j}
(C_\phi^{A'}(x,y)-C_\tau^{A'}(x,y)).$$
Therefore, we get
$$\begin{array}{l}
C_\phi^A(x,y)=C_\phi^{A'}(x,y)+C_\phi^{A'}(x,y)\sum\limits_{s\geq1}
x^{sa_k}y^s\sum\limits_{j=0}^s\left(C_\tau^{A'}(x,y)\right)^j\left(C_\nu^{A'}(x,y)\right)^{s-j}\\
\qquad\qquad\qquad\qquad\qquad-\sum\limits_{s\geq1}
x^{sa_k}y^s\sum\limits_{j=1}^s
\left(C_\tau^{A'}(x,y)\right)^j\left(C_\nu^{A'}(x,y)\right)^{s+1-j},\end{array}$$
or equivalently,
$$C_\phi^A(x,y)=(C_\phi^{A'}(x,y)-x^{a_k}yC_\tau^{A'}(x,y)C_\nu^{A'}(x,y))\sum_{s\geq0}
x^{sa_k}y^s\sum_{j=0}^s\left(C_\tau^{A'}(x,y)\right)^j\left(C_\nu^{A'}(x,y)\right)^{s-j}.$$
Hence, using the identity $\displaystyle\sum_{n\geq0}
x^n\sum_{j=0}^n p^jq^{n-j}=\frac{1}{(1-xp)(1-xq)}$ we get the
desired result  (1). Using similar arguments and replacing $a_1$
by $a_k$, we obtain (2).
\end{proof}

For certain shuffle patterns $\phi$ we can compute the generating function $C_\phi^A(x,y)$ explicitly, using the recursion  given in Theorem~\ref{thsb}.

\begin{example}
Let $A=\{a_1,\ldots,a_k\}$ be any ordered set of positive integers
and $\phi=1'\mn2\mn1''$ {\rm(}resp. $\psi=2'\mn1\mn2''${\rm)}. Here
$\tau=\nu=1$, so $C_\tau^A(x,y)=C_\nu^A(x,y)=1$ for any $A$, since
only the empty composition avoids $\tau$. Hence,
$$C_\phi^A(x,y)=\frac{1}{(1-x^{a_k}y)^2}(C_\phi^{A-\{a_k\}}(x,y)-x^{a_k}y).$$
Also, $C_\phi^{\{a_1\}}(x,y)=\frac{1}{1-x^{a_1}y}$
 as for any $m$, only the composition $\underbrace{a_1a_1\ldots
a_1}_{m\ times}$ avoids $\phi$ and therefore,
$$C_{1'\mn2\mn1''}^A(x,y) = \frac{1}{\prod_{a\in A}(1-x^ay)^2}
-\sum_{a\in A}\frac{x^ay}{\prod_{a\leq b\in A}(1-x^by)^2},$$ the result obtained in Equation ~(\ref{eqex1}) directly. Likewise, we obtain
$$C_{2'\mn1\mn2''}^A(x,y) = \frac{1}{\prod_{a\in A}(1-x^ay)^2}
-\sum_{a\in A}\frac{x^ay}{\prod_{a\geq b\in A}(1-x^by)^2}.$$
\end{example}

We now give two corollaries to Theorem~\ref{thsb}.

\begin{corollary}\label{thsc}
Let $\phi=\tau\mn \ell \mn\nu$ {\rm(}resp.
$\phi=\tau\mn1\mn\nu${\rm)} be a shuffle pattern, and let
$f(\phi)=f_1(\tau)\mn \ell \mn f_2(\nu)$ {\rm(}resp.
$f(\phi)=f_1(\tau)\mn1\mn f_2(\nu)${\rm)}, where
$f_1,f_2\in\{R,I\}$ are any trivial bijections. Then $\phi\equiv f(\phi)$.
\end{corollary}
\begin{proof}
Using Theorem~\ref{thsb}, and the fact that the number of
compositions in $C_{n;m}^A$ avoiding $\tau$ (resp. $\nu$) and
$f_1(\tau)$ (resp. $f_2(\nu)$) have the same generating functions, we get the desired
result.
\end{proof}

\begin{corollary}\label{thsd}
For any shuffle pattern $\tau\mn \ell \mn\nu$ {\rm(}resp.
$\tau\mn1\mn\nu${\rm)}, we have $\tau\mn \ell \mn\nu\equiv
\nu\mn\ell \mn\tau$ {\rm(}resp. $\tau\mn 1\mn\nu\equiv\nu\mn
1\mn\tau${\rm)}.
\end{corollary}
\begin{proof}
Corollary~\ref{thsc} yields that the pattern $\tau\mn \ell
\mn\nu$ (resp. $\tau\mn1\mn\nu$) is equivalent to the pattern
$\tau\mn \ell \mn R(\nu)$ (resp. $\tau\mn1\mn R(\nu)$), which is
equivalent to the pattern $R(\tau\mn \ell \mn R(\nu))=\nu\mn \ell
\mn R(\tau)$ (resp. $R(\tau\mn1\mn R(\nu))=\nu\mn1\mn R(\tau)$).
Finally, we use Corollary~\ref{thsc} one more time to get the
desired result.
\end{proof}

\subsection{Multi-patterns in compositions.}\label{sec4}

We now look at the second class of patterns. Recall that
a multi-pattern is of the form $\tau=\tau_1\mn\tau_2\mn\cdots\mn\tau_s$, where
$\{\tau_1,\dots,\tau_s\}$ is a set of consecutive patterns
and each letter of $\tau_i$ is incomparable with any letter of
$\tau_j$ whenever $i\neq j$.

The simplest non-trivial example of a multi-pattern is the pattern
$\phi=1\mn1'2'$. To avoid $\phi$ is the same as to avoid the
patterns $1\mn12$, $1\mn23$, $2\mn12$, $2\mn13$, and $3\mn12$
simultaneously. To count the number of compositions in
$C_{n;m}^A(1\mn1'2')$, we choose the leftmost letter of $\sigma \in
C_{n;m}^A(1\mn1'2')$ in $k$ ways, namely $a_1,\ldots,a_k$, and
observe that all the other letters of $\sigma$ must be in
non-increasing order. Hence,
$$C_{1\mn1'2'}^A(x,y)=1+\frac{y\sum_{a\in A}x^a}{\prod_{a\in A}(1-x^{a}y)}.$$

More generally, using
Lemma~\ref{praa} and Theorem~\ref{thmpc}, we get the following theorem that is the
basis for calculating the number of compositions that avoid a
multi-pattern, and therefore is the main result for
multi-patterns in this paper.

\begin{theorem}\label{coaa}
Let $A=\{a_1,\ldots,a_k\}$ be any ordered finite set of positive
integers and let $\tau=\tau_1\mn\tau_2\mn\cdots\mn\tau_s$ be a
multi-pattern. Then
$$C_\tau^A(x,y)=\sum_{j=1}^s C_{\tau_j}^A(x,y)\prod_{i=1}^{j-1}
\left[\left(y\sum_{a\in A}x^a-1\right)C_{\tau_i}^A(x,y)+1\right].$$
\end{theorem}

\begin{example}\label{coroldescents}
Let $A=\{a_1,\ldots,a_k\}$ be any ordered set of positive
integers. Let $\tau=\tau_1\mn\tau_2\mn\cdots\mn\tau_s$ be a
multi-pattern such that $\tau_j$ is equal to either $12$ or $21$,
for $j=1,2,\ldots,s$. It is easy to see that
$C_{12}^A(x,y)=C_{21}^A(x,y)=\frac{1}{\prod_{a\in A}(1-x^ay)}$ and we obtain from Theorem~\ref{coaa}
$$C_\tau^A(x,y)=\frac{1-\left(1+\frac{y\sum_{a\in A}x^a-1}{\prod_{a\in A}(1-x^ay)}\right)^s}
{1-y\sum_{a\in A}x^a}.$$
\end{example}

Using arguments similar to those in the proof of
\cite[Theorem~4.1]{KM} we get the following theorem which is an
analogue to~\cite[Theorem~21]{Ki2} and \cite[Theorem~4.1]{KM}.

\begin{theorem}\label{thmpa}
Let $\tau=\tau_0\mn\tau_1$ and $\phi=f_1(\tau_0)\mn f_2(\tau_1)$,
where $f_1$ and $f_2$ are any of the trivial bijections. Then
$\tau\equiv\phi$.
\end{theorem}

\begin{proof}
First we prove that  $\tau=\tau_0\mn\tau_1\equiv\tau_0\mn
f(\tau_1)$, where $f$ is a trivial bijection. Suppose that
$\sigma=\sigma_1\sigma_2\sigma_3$ avoids $\tau$ and
$\sigma_1\sigma_2$ has exactly one occurrence of $\tau_0$, namely
$\sigma_2$. Then $\sigma_3$ must avoid $\tau_1$, so $f(\sigma_3)$
avoids $f(\tau_1)$ and $\sigma_f=\sigma_1\sigma_2 f(\sigma_3)$
avoids $\phi$. The converse is also true, if $\sigma_f$ avoids
$\phi$ then $\sigma$ avoids $\tau$. Since any composition either avoids
$\tau_0$ or can be factored as above, we have a bijection between
the class of compositions avoiding $\tau$ and the class of compositions avoiding $\phi$.
Thus $\tau_0\mn\tau_1\equiv\tau_0\mn
f(\tau_1)$. Using this result as well as the properties of
trivial bijections we get
$$\begin{array}{l}
\tau\equiv\tau_0\mn f_2(\tau_1)\equiv R(\tau_0\mn
f_2(\tau_1))\equiv R(f_2(\tau_1))\mn R(\tau_0)\equiv\\
\qquad\qquad\qquad\qquad\qquad\equiv R(f_2(\tau_1))\mn
f_1(R(\tau_0))\equiv R(f_2(\tau_1))\mn R(f_1(\tau_0))\equiv
f_1(\tau_0)\mn f_2(\tau_1).
\end{array}$$
\end{proof}

\begin{corollary} \label{case2} The multi-patterns $\tau_1\mn\tau_2$ and $\tau_2\mn\tau_1$ are equivalent.
\end{corollary}

\begin{proof}
From Theorem~\ref{thmpa}, using the properties of the trivial bijection
$R$, we get
$$\tau_1\mn\tau_2\equiv \tau_1\mn R(\tau_2)\equiv R(R(\tau_2))\mn
R(\tau_1)\equiv \tau_2\mn R(R(\tau_1))\equiv \tau_2\mn\tau_1.$$
\end{proof}

We can obtain an even more general result.

\begin{theorem}\label{thmpb}
Suppose we have multi-patterns $\tau=\tau_1\mn\tau_2\mn\cdots\mn\tau_s$ and
$\phi=\phi_1\mn\phi_2\mn\cdots\mn\phi_s$, where $\tau_1\tau_2\ldots\tau_s$
is a permutation of $\phi_1\phi_2\ldots\phi_s$. Then $\tau\equiv\phi$.
\end{theorem}

\begin{proof} We use induction on $k$. For $s = 2$, the statement follows from Corollary~\ref{case2}. Suppose the statement is true for all $k < s$. If the composition $\sigma$ has no occurrences of $\phi_1$, then it obviously avoids both $\tau$ and $\phi$. Otherwise, we can write $\sigma = \sigma_1\sigma_2\sigma_3$, where $\sigma_1\sigma_2$ quasi-avoids $\phi_1$. Then $\sigma_3$ has to avoid $\phi_2\mn\cdots\mn\phi_s$. Since the $\phi_i$ are incomparable, it is irrelevant from which letters $\sigma_1\sigma_2$ is built, and we can apply the inductive hypothesis to $\phi_2\mn\cdots\mn\phi_s$. We can rearrange $\phi_2', \ldots \phi_k'$ of $\phi_2, \ldots \phi_k$ in such a way that the blocks in $\tau_1\tau_2\ldots\tau_s$ corresponding to $\phi_2,\ldots,\phi_s$ are arranged in the same order as the $\tau$'s. Then
\begin{eqnarray} \label{equivs} \phi=\phi_1\mn\phi_2\mn\cdots\mn\phi_s\equiv\phi_1\mn\phi_2'\mn\cdots\mn\phi_s'\equiv R(\phi_s')\mn\cdots\mn R(\phi_2')\mn R(\phi_1).
\end{eqnarray}
Now we consider two cases: Either $\tau_k \ne \phi_1$ or $\tau_k = \phi_1$. In the first case, we apply the hypothesis to the pattern $R(\phi_s')\mn\cdots\mn R(\phi_2')\mn R(\phi_1)$, with the role of $\phi_1$ played by $R(\phi_s')$. Thus, we can move the pattern $R(\phi_1)$ to the correct place somewhere to the left of $R(\phi_2')$, then apply the bijection $R$ to obtain that $\tau \equiv \phi$. In the second case, we obtain
$$\phi \equiv R(\phi_s')\mn\cdots\mn R(\phi_2')\mn R(\phi_1) \equiv R(\phi_s')\mn\cdots\mn R(\phi_1)\mn R(\phi_2')\equiv \phi_2'\mn \phi_1 \mn \cdots \mn \phi_s' \equiv  \phi_2'\mn \phi_s' \mn \cdots \mn \phi_1 = \tau.
$$
The first equivalence follows from ~(\ref{equivs}); the second one follows from the inductive hypothesis. Applying the bijection $R$ together with $R(R(x))=x$ and the inductive hypothesis once more gives the remaining equivalences.
\end{proof}

\section{Non-overlapping occurrences of POPs in
compositions}\label{non-overlap}

Kitaev~\cite{Ki2} and Mendes and Remmel \cite{Me,MR} proved the
following result on the distribution of non-overlapping patterns in
permutations: Let $\tau\mn$nlap($\sigma$) be the maximum number of
non-overlapping occurrences of a consecutive pattern $\tau$ in a
permutation $\sigma$ where two occurrences of $\tau$ are said to
overlap if they contain any of the same integers. Then
\begin{equation}\label{nov}\sum_{n=0}^{\infty}\frac{x^n}{n!}\sum_{\sigma\in
S_n}y^{\tau\mn\mbox{nlap}(\sigma)}=\frac{A(x)}{(1-y)+y(1-x)A(x)},\end{equation}
where $A(x)=\sum_{n=0}^{\infty}\frac{x^n}{n!}|\sigma\in S_n\ :\
\sigma\ \mbox{avoids}\ \tau|$. In other words, if the exponential
generating function for the number of permutations in $S_n$ avoiding
$\tau$ is known, then so is the bivariate generating function for
the entire distribution of~$\tau$-nlap. Kitaev and Mansour
\cite[Theorem 5.1]{KM} found an analogue to~(\ref{nov}) in case of
words. We now prove a corresponding result for compositions.

Let $\tau$ be an arbitrary consecutive pattern. We say that two
patterns {\em overlap} in a composition  if they contain any of the
same letters of the composition. Using Theorem~\ref{coaa} for the
multi-pattern $\tau \mn \tau \mn \cdots \mn \tau$ allows us to
obtain the generating function for the entire distribution of the
maximum number of non-overlapping occurrences of a pattern $\tau$ in
compositions.

The simplest consecutive pattern is a descent (or drop) in a
composition, which occurs at position $i$ if
$\sigma_i>\sigma_{i+1}$. Clearly, two descents $i$ and $j$  overlap
if $j=i+1$. In particular, we can define the statistic {\em maximum
number of non-overlapping descents}, or $\MND$, in a composition.
For example, $\MND(333211)=1$ whereas $\MND(13321111432111)=3$
(namely 32, 43 and 21). Obviously, this statistic, maximum number of
non-overlapping patterns, can be defined for any consecutive pattern
$\tau$, and we obtain the following result.

\begin{theorem}\label{dstmnd}
Let $A$ be any ordered set of positive integers and let $\tau$ be a
consecutive pattern. Then
$$\sum_{n,m\geq0}\sum_{\sigma\in C_{n;m}^A} t^{\tau\mn\mbox{nlap}(\sigma)}x^ny^m=
\frac{C_\tau^A(x,y)}{1-t\left[\left(y\sum_{a\in A}x^a
-1\right)C_\tau^A(x,y)+1\right]},$$ where $\tau\mn$nlap($\sigma$) is
the maximum number of non-overlapping occurrences of $\tau$ in
$\sigma$.
\end{theorem}
\begin{proof}
We fix a natural number $s$ and consider the multi-pattern
$\Phi_s=\tau\mn\tau\mn\cdots\mn\tau$ with $s$ copies of $\tau$. If a
composition avoids $\Phi_s$ then it has at most $s-1$
non-overlapping occurrences of $\tau$. Theorem~\ref{coaa} yields
$$C_{\Phi_s}^A(x,y)=\sum_{j=1}^s C_{\tau}^A(x,y)\prod_{i=1}^{j-1} \left[\left(y\sum_{a\in A}x^a-1\right)C_{\tau}^A(x,y)+1\right].$$
Therefore, the generating function for the number of compositions that have exactly $s$
non-overlapping occurrences of the pattern $\tau$ is given by
$$C_{\Phi_{s+1}}^A(x,y)-C_{\Phi_s}^A(x,y)=C_\tau^A(x,y)\left[\left(y\sum_{a\in A}x^a-1\right)C_\tau^A(x,y)+1\right]^s.$$
Hence,
$$\sum_{n,m\geq0}\sum_{\sigma\in C_{n;m}^A}t^{\tau\mn\mbox{nlap}(\sigma)}x^ny^m=\sum_{s\geq0}
t^sC_\tau^A(x,y)\left[\left(y\sum_{a\in A}x^a
-1\right)C_\tau^A(x,y)+1\right]^s,$$ or, equivalently,
$$\sum_{n,m\geq0}\sum_{\sigma\in C_{n;m}^A}t^{\tau\mn\mbox{nlap}(\sigma)}x^ny^m=\frac{C_\tau^A(x,y)}{1-t\left[\left(y\sum_{a\in A}x^a
-1\right)C_\tau^A(x,y)+1\right]}.$$
\end{proof}

Note that Theorem~\ref{dstmnd} is a $q$-analogue to \cite[Theorem
5.1]{KM}, which is the main result of~\cite{KM} (set $x=1$ to get the
result for words). We use Theorem~\ref{dstmnd} to
obtain the distribution for $M\!N\!D$, the maximum number of
non-overlapping descents.

\begin{example}\label{ex5.2}
Let $A$ be any ordered set of positive integers. If we consider
descents {\rm(}the pattern $12${\rm)} then
$C_{12}^A(x,y)=\frac{1}{\prod_{a\in A}(1-x^ay)}$, hence the
distribution of $M\!N\!D$ is given by the formula:
$$\sum_{n,m\geq0}\sum_{\sigma\in
C_{n;m}^A}t^{12\mn\mbox{nlap}(\sigma)}x^ny^m =\frac{1}{\prod_{a\in
A} (1-x^ay)+t\left(1-y\sum_{a\in A}x^a-\prod_{a\in
A}(1-x^ay)\right)}.$$ Specifically, the distribution of $M\!N\!D$ on
the set of compositions of $n$ with parts in $A=\{1,2\}$ is given by
$$\frac{1}{(1-x)(1-x^2)-x^3t}=\sum_{s\geq0}\frac{x^{3s}}{(1-x)^{2s+2}(1+x)^{s+1}}t^s.$$
\end{example}


\end{document}